\documentclass[12pt]{amsart}
\usepackage{amssymb}
\usepackage{graphicx}
\newtheorem{fact}{Fact}[section]
\newtheorem{lemma}[fact]{Lemma}
\newtheorem{theorem}[fact]{Theorem}
\newtheorem{corollary}[fact]{Corollary}
\theoremstyle{definition}
\newtheorem{definition}[fact]{Definition}

\newcommand{\la}{\lambda}
\DeclareMathOperator{\maj}{maj}
\DeclareMathOperator{\Inv}{Inv}
\DeclareMathOperator{\Des}{Des}
\DeclareMathOperator{\inv}{inv}

\title[An Extension of the Foata Map to Standard Young Tableaux] {An Extension of the Foata Map to Standard Young Tableaux}
\author[J. Haglund and L. Stevens]
{J. Haglund$^{\,\star,1}$ and L. Stevens$^{\,\diamond}$}
\thanks{$^1$ Supported in part by NSA grant MSPF-02G-193 and NSF grant DMS-0553619}
\begin{document}
\maketitle
\begin{center}
{{\it $^\star$ Department of Mathematics, University of Pennsylvania, Philadelphia, PA 19104-6395}, \texttt{jhaglund@math.upenn.edu}}
\medskip

{\it $^\diamond$ Department of Mathematics, University of California, San Diego, La Jolla, CA 92093-0112}, \texttt{stevens@math.ucsd.edu}
\end{center}
\centerline{October, 2006}
\begin{abstract} We define an inversion statistic on standard Young tableaux.  We prove that this statistic has the same distribution over $SYT(\la)$ as the major index statistic by exhibiting a bijection on $SYT(\la)$ in the spirit of the Foata map on permutations.
\end{abstract}
\section{Introduction}
A permutation statistic is a combinatorial rule which associates an element of $\mathbb{N}$ to each element of the symmetric group $S_{n}$.  Let $\sigma=\sigma_{1}\sigma_{2}\cdots\sigma_{n}$ denote the element of $S_{n}$ which sends $i$ to $\sigma_{i}$ for $1\leq i\leq n$.  An inversion of $\sigma$ is a pair $(i,j),\,1\leq i<j\leq n$ such that $\sigma_{i}>\sigma_{j}$.  A descent of $\sigma$ is an integer $i,\,1\leq i\leq n-1$, for which $\sigma_{i}>\sigma_{i+1}$.  The inversion statistic $\inv(\sigma)$ is defined to be the number of inversions of $\sigma$ and the major index statistic $\maj(\sigma)$ is defined to be the sum of the descents of $\sigma$, i.e.
$$\inv(\sigma)=\sum_{i<j,\,\sigma_{i}>\sigma_{j}}1,\qquad\maj(\sigma)=\sum_{i,\,\sigma_{i}>\sigma_{i+1}}i.$$
Major P. MacMahon \cite{MacM} introduced the major index statistic and proved that, remarkably, its distribution over $S_{n}$  is equal to the distribution of the inversion statistic over $S_n$.  
This raised the question of constructing a canonical bijection $\phi:S_{n}\to S_{n}$ such that $\maj(\sigma)=\inv(\phi(\sigma))$.  Foata \cite{Foa68} found such a map. (See also
\cite{FoSc78}).  Note that the previous two statistics can be defined on each class $\tilde{C}$ of permutations of a sequence (with repetitions) $1^{m_1}2^{m_2}\cdots r^{m_r}$.  The equidistribution of $\inv$ and $\maj$ also holds for such a class $\tilde{C}$, a result proved by MacMahon \cite{MacM} and reproved by means of the Foata map \cite{Foa68}, when adequately defined on $\tilde{C}$.  In this paper, we only need the case of permutations without repetitions.

There is a version of $\text{maj}$ for tableau which plays a
prominent role in symmetric function theory (see for example \cite[Chapter 7]{Sta} or
\cite[Chapter 1]{Mac}).   One could also ask whether there is a natural version of
$\text{inv}$ which would play a similar role.  In this article we advance a candidate 
tableau inversion statistic, which is defined in terms of a generalization of
Foata's map.  The statistic was discovered in the course of studying Macdonald 
polynomials (see the remark at the end of Section \ref{bijection}).

In section one, we review the algorithm describing the Foata map and the necessary background on tableaux.  In section two we define the inversion statistic $\Inv$ and maps on tableaux that will feature in our Foata-type map.  In section three we introduce our Foata-type map and prove that it is a bijection which sends a standard Young tableau with a given $\Inv$ to a standard Young tableau with the same major index. Section \ref{extensions}
contains some remarks about extending our map to skew shapes, and 
in what sense it
generalizes Foata's map.

\section{Definitions and Background}
We begin by reviewing Foata's map $\phi:S_{n}\to S_{n}.$  It can be described as follows.  If $n\leq2,$ $\phi(\sigma)=\sigma$.  If $n>2$, we add numbers to $\phi$ one at a time: begin by setting $\phi^{(1)}=\sigma_{1}$, $\phi^{(2)}=\sigma_{1}\sigma_{2}$.   To find $\phi^{(3)}$, start with $\sigma_{1}\sigma_{2}\sigma_{3}$.  Then if $\sigma_{3}>\sigma_{2}$, draw a bar after each element of $\sigma_{1}\sigma_{2}\sigma_{3}$ which is less than $\sigma_{3}$, while if  $\sigma_{3}<\sigma_{2}$, draw a bar after each element of $\sigma_{1}\sigma_{2}\sigma_{3}$ which is greater than $\sigma_{3}$.  Also add a bar before $\sigma_{1}$.  For example, if $\sigma=4137562$, we now have $\vert41\vert3.$  Now regard the numbers between two consecutive bars as ``blocks", and in each block, move the last element to the beginning, and finally remove all of the bars.  We end up with $\phi^{(3)}=143.$

Proceeding inductively, we begin by adding $\sigma_{i}$ to the end of $\phi^{(i-1)}$.  Then if $\sigma_{i}>\sigma_{i-1}$, draw a bar after each element of $\phi^{(i-1)}$ which is less than $\sigma_{i}$, while if $\sigma_{i}<\sigma_{i-1}$, draw a bar after each element of $\phi^{(i-1)}$ which is greater than $\sigma_{i}$.  Also draw a bar before $\phi_{1}^{(i-1)}$.  Then in each block, move the last element to the beginning, and finally remove all of the bars.  If $\sigma=4137562$, the successive stages of the algorithm yield
\begin{align*}
143&=\phi^{(3)}\\
\vert1\vert4\vert3\vert7\mapsto1437&=\phi^{(4)}\\
\vert1437\vert5\mapsto71435&=\phi^{(5)}\\
\vert71\vert4\vert3\vert5\vert6\ \mapsto174356&=\phi^{(6)}\\
\vert17\vert4\vert3\vert5\vert6\vert2\mapsto7143562&=\phi^{(7)},&
\end{align*}
so $\phi(4137562)=7143562$.  Note that $\maj(4137562)=11=\inv(7143562)$.  
\begin{theorem}\cite{Foa68}
The map  $\phi:S_{n}\to S_{n}$ is a bijection and for all $\sigma\in S_{n}$, $\maj(\sigma)=\inv(\phi(\sigma))$.  
\end{theorem}
The inverse to Foata's map is as follows.  Let $\psi=\phi^{-1}(\sigma)$.  If $\sigma_{n}>\sigma_{1}$, draw a bar before each number in $\sigma$ which is less than $\sigma_{n}$ and also before $\sigma_{n}$.  If $\sigma_{n}<\sigma_{1}$, draw a bar before each number in $\sigma$ which is greater than $\sigma_{n}$ and also before $\sigma_{n}$.  Next move each number at the beginning of a block to the end of the block.  Remove the bars to obtain $\psi^{(1)}$.
The last letter of $\psi$ is now fixed.  Now compare $\psi^{(1)}_{n-1}$ with $\psi^{(1)}_{1}$.  Create blocks as above, drawing a bar before $\psi^{(1)}_{n-1}$.  Proceed in this manner.  For example, if   
$\sigma=7143562$, the successive stages of the algorithm yield
\begin{align*}
\vert71\vert4\vert3\vert5\vert6\vert2\mapsto1743562&=\psi^{(1)}\\
\vert17\vert4\vert3\vert5\vert62\mapsto7143562&=\psi^{(2)}\\
\vert7143\vert562\mapsto1437562&=\psi^{(3)}\\
\vert1\vert4\vert3\vert7562\mapsto1437562&=\psi^{(4)}\\
\vert14\vert37562\mapsto4137562&=\psi^{(5)}\\
4137562=\psi^{(6)}&=\psi^{(7)},
\end{align*}
so $\phi^{-1}(7143562)=4137562$.  It turns out that the generalization of $\phi^{-1}$ to tableaux is simpler to describe than the generalization of $\phi$.

Next we review the necessary background on tableaux.
A partition of a positive integer $n$ is a sequence $\la=(\la_{1},\dots,\la_{k})$ of positive integers such that $\sum_{i=1}^{k} \la_{i}=n$ and $\la_{1}\geq\cdots\geq\la_{k}$.  If $\la$ is a partition of $n$ we write $\la\vdash n$ or $|\la|=n$.  The $\la_{i}$ are called the parts of $\la$.

We represent $\la=(\la_{1},\dots,\la_{k})$ pictorially as follows. 
We assign (row,column)-coordinates to unit squares in the first quadrant, 
obtained by permuting the $(x,y)$ coordinates of the upper right-hand corner of the square, so the lower left-hand square has coordinates $(1,1)$, the square above it $(2,1)$, etc., and a square $(i,j)$ has ordinate $i$ and abscissa $j$.
Then the Ferrer's diagram of $\la$, which we also denote by $\la$, is the set of squares, or ``cells"  $\{(i,j)\in\mathbb{N}\times\mathbb{N}:1\leq i\leq k,\,1\leq j\leq\la_{i}\}$.  That is, the diagram of $\la$ consists of $k$ left-justified rows of squares in the first quadrant of the $xy$-plane with $\la_{i}$ squares in the $i$th row from the bottom.  See Figure $1$.

\begin{figure}\label{ferrer}
\begin{center}
\includegraphics{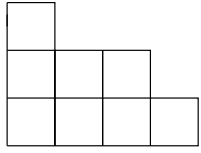}\end{center}
\caption{The Ferrer's diagram of the partition $(4,3,1)$}
\end{figure}

The conjugate partition of $\la=(\la_{1},\dots,\la_{k})$ is $\la'=(\la_{1}',\dots,\la_{j}')$, where $\la_{i}'$ is the number of parts of $\la$ which are greater than or equal to $i$.  Geometrically, the diagram of $\la'$ is the reflection of the diagram of $\la$ about the line $y=x$.

If $\la\vdash n$, a standard Young tableau of shape $\la$ is a filling of the diagram of $\la$ with the numbers $1,\dots,n$ such that in each row, the numbers are increasing from left to right, and in each column, the numbers are increasing from bottom to top.  The number filling a cell is called the content of the cell.  We let $c(i,j)$ denote the content of the cell $(i,j)$.

We let $SYT(\la)$ denote the set of standard Young tableaux of shape $\la$.  
For $T\in SYT(\la)$, the conjugate tableau $T'\in SYT(\la')$ is obtained by filling the cell $(i,j)$ with the content of the cell $(j,i)$ of $T$.

\section{An $\Inv$ Statistic on Standard Young Tableaux}
For a standard Young tableau $T$, the major index of $T$ is given by
$$\maj(T)=\sum_{i\in \Des(T)}i,$$
where $\Des(T)=\{i\,|\,i+1 \text{ is in a row above the row containing $i$ in $T$}\}$. 
Here we define another statistic on standard Young tableaux which we call the $\Inv$ statistic.  Let $(i,j),(k,l)\in T$ such that $c(i,j)>c(k,l)$.  The cells $(i,j)$ and $(k,l)$ form an ``inversion pair"  if $(k,l)$ is weakly SE of $(i,j)$, and they do not form an inversion pair if $(k,l)$ is weakly NW of $(i,j)$.  If $(k,l)$ is strictly SW of $(i,j)$, then whether or not they form an inversion pair is determined by a path in $T$ called the ``inversion path" of $(i,j)$.  By strictly SW, we mean Southwest but not due South or due West, i.e. $k<i$ and $l<j$.  Below, we construct the set of of inversion paths for $T$.  In the construction we define maps on standard Young tableaux which will be used in our extension of the Foata map.         

\begin{definition} Let $T\in SYT(\la)$, where $\la\vdash n$.  For $k=1,\dots,n$, let $\pi(T,k)$ be the path in $T$ constructed according to the following algorithm:  

Start at the lower left-hand corner of the cell with content $k$.  Suppose this cell has coordinates $(i,j)$.  If $i=1$ or $j=1$, proceed in a straight line until you reach the origin.  If not, compare the contents of the cells $(i-1,j)$ and $(i,j-1)$.  If $c(i-1,j)>c(i,j-1)$, take one unit step South.  If $c(i-1,j)<c(i,j-1)$, take one unit step West.  Arriving at the lower left-hand corner of a new cell $(i',j')$, iterate the algorithm.  Proceed until you reach the origin.
\end{definition} 
For example, let $T$ be the SYT in Figure $2$ (the fact that the shape is a square is
coincidental).  To form the inversion path $\pi (T,16)$, we start at the lower left-hand
corner of the square containing $16$, and draw a unit line segment left, since the square to the left of $16$ contains $15$, which is greater than the number in the square below
(namely $14$).  Next we compare the numbers in the squares to the left and below the square
containing $15$.  The largest is $13$, so we now extend our inversion path
by drawing another unit line segment, in this case downwards, ending at the lower left-hand corner of the square containing $13$.  We continue comparing squares left and below, moving one unit in the direction of the
largest, until we reach the edge of the tableau (the lower left-hand corner of the square
containing $3$), at which time we go in the only direction we can, namely straight down,
until we reach a square with no left or bottom neighbors (the square containing $1$).

\begin{figure}\label{path}
\begin{center}
{\includegraphics{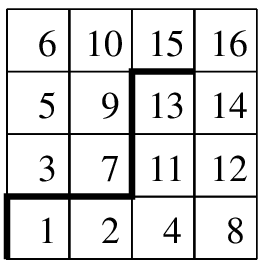}}\end{center}
\caption{A standard Young tableau $T$ and the path $\pi(T,16)$}
\end{figure}

Note that the path $\pi(T,k)$ partitions the cells of T with content less than $k$ into two sets: cells which are weakly NW of $\pi(T,k)$ and cells which are weakly SE of $\pi(T,k)$.  
We will say that cells weakly NW of $\pi(T,k)$ are above $\pi(T,k)$ and cells weakly SE of $\pi(T,k)$ are below $\pi(T,k)$.  
Two cells are on the same side of $\pi(T,k)$ if they are both above or both below $\pi(T,k)$.  Note also that since $\pi(T,k)$ is constructed taking only South and West steps from the cell with content $k$, all cells weakly SE of the cell with content $k$ are always below $\pi(T,k)$ and all cells weakly NW of the cell with content $k$ are always above $\pi(T,k)$.
\begin{definition} Let $T\in SYT(\la)$, where $\la\vdash n$.  For $k=1, \dots, n$, let $\Psi_{k}(T)$ be the tableau constructed according to the following algorithm:    

Construct the path $\pi(T,k)$.  Partition the cells of $T$ with content less than $k$ into ``blocks" $\{(i_{1},j_{1}),\dots,(i_{m},j_{m})\}$ of maximal length such that 
\begin{enumerate}
\renewcommand{\theenumi}{\roman{enumi}}\renewcommand{\labelenumi}{(\theenumi)}
\item $c(i_l,j_l)=c(i_{l+1},j_{l+1})-1$ for $1\leq l\leq m-1$. 
\item$(i_{1},j_{1})$ is on the same side of $\pi(T,k)$ as $(1,1)$,  and $(i_{l},j_
{l})$ is on the other side of $\pi(T,k)$ for $2\leq l\leq m$,    
\end{enumerate}
Note that blocks may consist of only one cell.  Next perform the following ``cycling" procedure on the contents of the cells in each block:  
\begin{enumerate}
\renewcommand{\theenumi}{\roman{enumi}}\renewcommand{\labelenumi}{(\theenumi)}
\item Replace $c(i_{1},j_{1})$ by $c(i_{m},j_{m})$.   
\item For $2\leq l\leq m$, replace $c(i_{l},j_{l})$ by  $c(i_{l},j_{l})-1$.  \end{enumerate}
The resulting tableau is $\Psi_{k}(T)$.   
\end{definition} 
In the next section, we prove that $\Psi_{k}(T)\in SYT(\la)$.  Note that the maps $\Psi_{2},\,\Psi_{1}$ are the identity. 

Simply put, the blocks for the path $\pi(T,k)$ consist of cells containing consecutive numbers such that the smallest number is in a cell on the same side of the path as the cell containing $1$ and all of the rest of the cells are on the other side of the path.  In the example in Figure $2$, if we denote cells by their contents, then the blocks for $\pi(T,16)$ are $\{1\},\,\{2,3\},\,\{4,5,6,7\},\,\{8,9,10\},\,\{11\},\,\{12\},\,\{13\},\,\{14,15\}$.  

If the contents of the cells in a block are $a,a+1,\dots,a+j$, then the cycling procedure sends $a$ to the cell which had content $a+1$, etc. and finally $a+j$ to the cell which had content $a$.  See Figure $3$.  A key observation about the cycling procedure is that it preserves the relative order of the contents of cells in different blocks.  

\begin{figure}\label{map}
\begin{center}
{\includegraphics{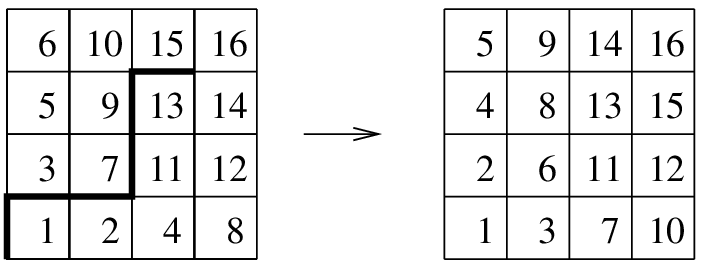}}\end{center}
\caption{The tableau on the right is $\Psi_{16}(T)$}
\end{figure}
\begin{definition}
The collection of paths 
$$\{\pi(T,n),\pi(\Psi_{n}(T),n-1),\dots,\pi(\Psi_{3}\circ\cdots\circ\Psi_{n}(T),2)\}$$ 
is the  set of inversion paths for $T$.  For $(i,j)\in T,\,(i,j)\neq(1,1)$, the inversion path of $(i,j)$ is the unique element of the set of inversion paths for $T$ starting at the lower left-hand corner of the cell $(i,j)$.  
\end{definition}
In our running example, the path $\pi(\Psi_{16}(T),15)$ is the inversion path for the cell $(3,4)$ (the cell with content $14$) in $T$.  See Figure $4$.  See also Figures $6$ and $7$ for a complete example.

\begin{figure}\label{path2}
\begin{center}
{\includegraphics{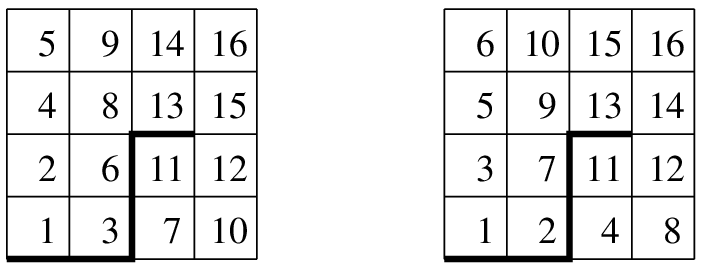}}\end{center}
\caption{ On the left: $\Psi_{16}(T)$ and the path $\pi(\Psi_{16}(T),15)$.  On the right: $T$ and the inversion path of the cell $(3,4)$.}
\end{figure}

Now we are ready to introduce the $\Inv$ statistic.
\begin{definition}
Let $T\in SYT(\la)$.  An ordered pair of cells $((i,j),(k,l))$ in $T$ is an inversion pair if $c(k,l)<c(i,j)$ and $(k,l)$ is below the inversion path of $(i,j)$.  $\Inv(T)$ is the total number of inversion pairs in $T$.
\end{definition} 
In our running example, the cell with content $16$ in $T$ forms inversion pairs with the cells with contents $1,2,4,8,11,12,13,14$ and the cell with content $14$ in $T$ forms inversion pairs with the cells with contents $4,8,11,12$.

\section{The bijection}
\label{bijection}
Let $\la\vdash n$, and let Let $\Psi=\Psi_{3}\circ\cdot\circ\Psi_{n}$.  One of our main results is:
\begin{theorem}\label{main} The map $\Psi:SYT(\la)\to SYT(\la)$ is a bijection and for $T\in SYT(\la)$, $\Inv(T)=\maj(\Psi(T))$.
\end{theorem}
Most of this section is devoted to the proof of Theorem \ref{main}.  First we prove that $\Psi_{k}:SYT(\la)\to SYT(\la)$ .  We will use the following lemma, which helps us understand the geometry of the blocks for a given path.    
\begin{lemma}\label{geom}   Let $T\in SYT(\la)$.  Suppose that the cells with contents $a,a+1,\dots,a+m$ form a block for $\pi(T,k)$.  If the cell $(1,1)$ is below $\pi(T,k)$, then the cells containing $a+1,\dots,a+m$ are strictly NW of the cell containing $a$.  If the cell $(1,1)$ is above $\pi(T,k)$, then the cells containing $a+1,\dots,a+m$ are strictly SE of the cell containing $a$. 
\end{lemma}
\begin{proof}  We assume $k$ is not in the first row or first column of $T$ because if it were, the lemma would be trivial as all cells would be on the same side of $\pi(T,k)$.  Suppose the cell $(1,1)$ is below $\pi(T,k)$.  Let $(i,j)$ be the cell with content $a$ and let $(r,s)$ be the cell with content $k$.  Since $(i,j)$ is below $\pi(T,k)$, $(i,j)$ must be either weakly SE or strictly SW of $(r,s)$.  In the first case, the remainder of the cells in the block must be in columns to the left of column $j$, since if they were not, they would be below $\pi(T,k)$.  These cells must also be in rows above row $i$ because all cells weakly SW of $(i,j)$ have content $<a$.  Next we consider the case in which $(i,j)$ is strictly SW of $(r,s)$.  
Since $\pi(T,k)$ consists of only South and West steps from $(r,s)$, every cell weakly SE of $(i,j)$ is also below $\pi(T,k)$.  Thus, the remainder of the cells in the block must be in rows above row $i$.  As $c(i,j+1)>a$ and the numbers $a+1,\dots,a+m$ occur in rows above row $i$, $c(i,j+1)>a+m$.  Hence, every cell weakly NE of $(i,j+1)$ has content $>a+m$.  So the remainder of the cells in the block must be in column $j$ or in columns to left of column $j$.  Suppose $c(i+p,j)=a+l$, where $a+l$ is the smallest of the numbers $a+1,\dots,a+m$ in column $j$.  As $(i+p,j)$ is above $\pi(T,k)$ and $(i,j)$ is below $\pi(T,k)$, there must be a West step of the path across column $j$.  To take this West step, we must have had $c(i+q,j)>c(i+q-1,j+1)$ for some $q$ between $1$ and $p$.  But this is a contradiction since $a<c(i+q,j)<a+m$ and $c(i+q-1,j+1)>a+m$.  Thus the remainder of the cells in the block must be in columns to the left of column $j$.    

The second statement follows since $(1,1)$ is below $\pi(T,k)$ if and only if $(1,1)$ is above $\pi(T',k)$, and the map $T\mapsto T'$ sends cells strictly NW of a given cell in $T$ to cells strictly SE of the image of that cell in $T'$.
\end{proof}
\begin{theorem}
For $k=3,\dots,n$, if $T\in SYT(\la)$, then $\Psi_{k}(T)\in SYT(\la)$.
\end{theorem}
\begin{proof}
Since $\Psi_{k}$ acts as the identity on all cells with content greater than $k-1$, it suffices to prove the theorem for $k=n$.
As our argument involves the content of both cells in $T$ and in $\Psi_{n}(T)$, it will be convenient to use the notation $c(i,j)_{S}$ to denote the content of the cell $(i,j)$ in the tableau $S$.  Let $(i,j)\in T$.  Let $c(i,j)_{T}=p$ and $c(i,j)_{\Psi_{n}(T)}=p'$. 
We need to show that if $(i+1,j)\in T$, then $p'<c(i+1,j)_{\Psi_{n}(T)}$ and if $(i,j+1)\in T$, then $p'<c(i,j+1)_{\Psi_{n}(T)}$.  Assume $(i+1,j)\in T$.  Let $c(i+1,j)_{T}=r$ and $c(i+1,j)_{\Psi_{n}(T)}=r'$.  
If $(i,j)$ and $(i+1,j)$ were in different blocks for $\pi(T,n)$, then $p'<r'$ since $p<r$ and the cycling procedure preserves the relative order of the contents of cells in different blocks.  If $(i,j)$ and $(i+1,j)$ were in the same block, then by Lemma \ref{geom}, they had to be on the same side of $\pi(T,n)$.  But then $p'=p-1$ and $r'=r-1$.  So $p<r$ implies $p'<r'$.  A similar argument works for the case $(i,j+1)\in T$. 
\end{proof}
Next we show that maps $\Psi_{k}$ are bijections by explicitly constructing the inverse maps.  
Our starting point is the following lemma, which is key in reconstructing the path $\pi(\Psi_{k}^{-1}(S),k)$.

\begin{lemma}\label{maj} Let $T\in SYT(\la)$, where $\la\vdash n$.  For $k=2,\dots,n$, $(1,1)$ is below $\pi(T,k)$ if and only if $k-1\in \Des(\Psi_{k}(T))$.  
\end{lemma}
\begin{proof}
We may assume $k$ is not in the first row or column of $T$, since in these cases the lemma is trivial. 
Suppose $(1,1)$ is below $\pi(T,k)$.  Let $a$ be the largest number below $\pi(T,k)$.  So $a$ is in a row below $k$.  If $a=k-1$, then the cell containing $a$ is in a block by itself, and the cycling procedure does not change the contents of this cell.  If $a<k-1$, then the cells containing $a,\dots,k-1$ form a block, and the cycling procedure sends $k-1$ to the cell which contained $a$.  So in either case, $k-1$ is in a row below $k$ in  $\Psi_{k}(T)$.  Using a similar argument, one shows that if $(1,1)$ is above $\pi(T,k)$, then $k-1$ is not in a row below $k$ in $\Psi_{k}(T)$.  
\end{proof}

Now we are ready to describe the maps $\Phi_{k}=\Psi_{k}^{-1}$.  Let $S\in SYT(\la)$, and suppose $c(r,s)=k$ and $c(x,y)=k-1$.  We give an algorithm to obtain a new tableau from $S$ and a path $\tilde{\pi}(S,k)$ starting at the lower left-hand corner of the cell $(r,s)$ and ending at the origin:\\     

\noindent{\it Step 1}:
If $r=1$ or $s=1$, draw a straight line to the origin and stop. Otherwise:\\
If $k-1\in \Des(S)$: Find all ``simple blocks" $\{(i_{1},j_{1}),\dots,(i_{m},j_{m})\}$ of cells such that
\begin{enumerate}
\renewcommand{\theenumi}{\roman{enumi}}\renewcommand{\labelenumi}{(\theenumi)}
\item $c(r-1,s-1)<c(i_{l},j_{l})<k$ and $c(i_{l},j_{l})=c(i_{l+1},j_{l+1})+1$ for all $1\leq k\leq m$,
\item $(i_{1},j_{1})$ is weakly SE of $(r,s)$ and $(i_{2},j_{2}),\dots,(i_{m},j_{m})$ are weakly NW of $(r,s)$,
\item the cell with content $c(i_{m},j_{m})-1$ is weakly SE of $(r,s)$.  
\end{enumerate}
If $k-1\notin \Des(S)$:
Find all simple blocks $\{(i_{1},j_{1}),\dots,(i_{m},j_{m})\}$ of cells such that
\begin{enumerate}
\renewcommand{\theenumi}{\roman{enumi}}\renewcommand{\labelenumi}{(\theenumi)}
\item $c(r-1,s-1)<c(i_{l},j_{l})<k$ and $c(i_{l},j_{l})=c(i_{l+1},j_{l+1})+1$ for all $1\leq k\leq m$,
\item $(i_{1},j_{1})$ is weakly NW of $(r,s)$ and $(i_{2},j_{2}),\dots,(i_{m},j_{m})$ are weakly SE of $(r,s)$,
\item the cell with content $c(i_{m},j_{m})-1$ is weakly NW of $(r,s)$.  
\end{enumerate}

\noindent{\it Step 2:}
Perform the following ``reverse cycling" procedure on the contents of the cells of each block from Step $1$:
Replace $c(i_{1},j_{1})$ by $c(i_{m},j_{m})$, and for $2\leq l\leq m$, replace $c(i_{l},j_{l})$ by $c(i_{l},j_{l})+1$.

\noindent{\it Comment:}  In Step 1, we are looking for the blocks of $\pi(\Psi_{k}^{-1}(S),k)$  which are easy to identify, in particular those which are independent of the path.  In Step 2, we undo the cycling procedure of $\Psi_{k}$ on these blocks.  To identify the rest of the blocks, we must reconstruct the path.

\noindent{\it Step 3:} Construct a unit segment of $\tilde\pi(S,k)$ as follows:\\
If $k-1\in \Des(S)$: Go South if $(r-1,s)$ is in one of the blocks from Step $2$ and $c(r-1,s)>c(r,s-1)$.  Otherwise go West.\\
If $k-1\notin \Des(S)$: Go West if $(r,s-1)$ is in one of the blocks from Step $2$ and $c(r,s-1)>c(r-1,s)$.  Otherwise go South. 

\noindent{\it Step 4:} A portion of $\tilde{\pi}(S,n)$ has been constructed and ends at the lower left-hand corner of the cell $(u,v)$.  If $u=1$ or $v=1$, draw a straight line to the origin and stop.  Otherwise:\\
Find all simple blocks $\{(i_{1},j_{1}),\dots,(i_{m},j_{m})\}$ of cells such that
\begin{enumerate}
\renewcommand{\theenumi}{\roman{enumi}}\renewcommand{\labelenumi}{(\theenumi)}
\item for all $1\leq k\leq m$, $(i_{l},j_{l})$ was not in a simple block at any previous stage in the algorithm,
\item $c(u-1,v-1)<c(i_{l},j_{l})<k$ and $c(i_{l},j_{l})=c(i_{l+1},j_{l+1})+1$ for all $1\leq k\leq m$,
\item $(i_{1},j_{1})$ and $(x,y)$ are on the same side of the portion of the path constructed so far and $(i_{2},j_{2}),\dots,(i_{m},j_{m})$ are on the other side, 
\item the cell with content $c(i_{m},j_{m})-1$ and $(i_{1},j_{1})$ are on the same side of portion of the path constructed so far.  
\end{enumerate}
\noindent{\it Step 5:} Reverse cycle the contents of the cells of each block from Step $4$.\\
\noindent{\it Step 6:} Construct a unit segment of $\tilde\pi(S,k)$ as follows:\\
If $n-1\in \Des(S)$: Go South if  $(u-1,v)$ has already appeared in a simple block and $c(u-1,v)>c(u,v-1)$.  Otherwise go West.\\
If $n-1\notin \Des(S)$: Go West if $(u,v-1)$ has already appeared in a simple block and $c(u,v-1)>c(u-1,v)$.  Otherwise go South. \\
Arriving at the lower left-hand corner of a new cell, iterate the algorithm from Step $4$.\\

Observe that if a cell with content $a$ appears in a simple block, then all cells with content greater than $a$ and less than $n$ must have already appeared in simple blocks.  
Moreover, upon completion of the algorithm, every cell with content greater than $1$ and less than $n$ will appear in exactly one simple block.  

In Figure $5$, we give an example of the algorithm with $k=16$ applied to the tableau $\Psi_{16}(T)$ from Figure $3$.  
\begin{figure}\label{inverse}
\begin{center}
{\includegraphics{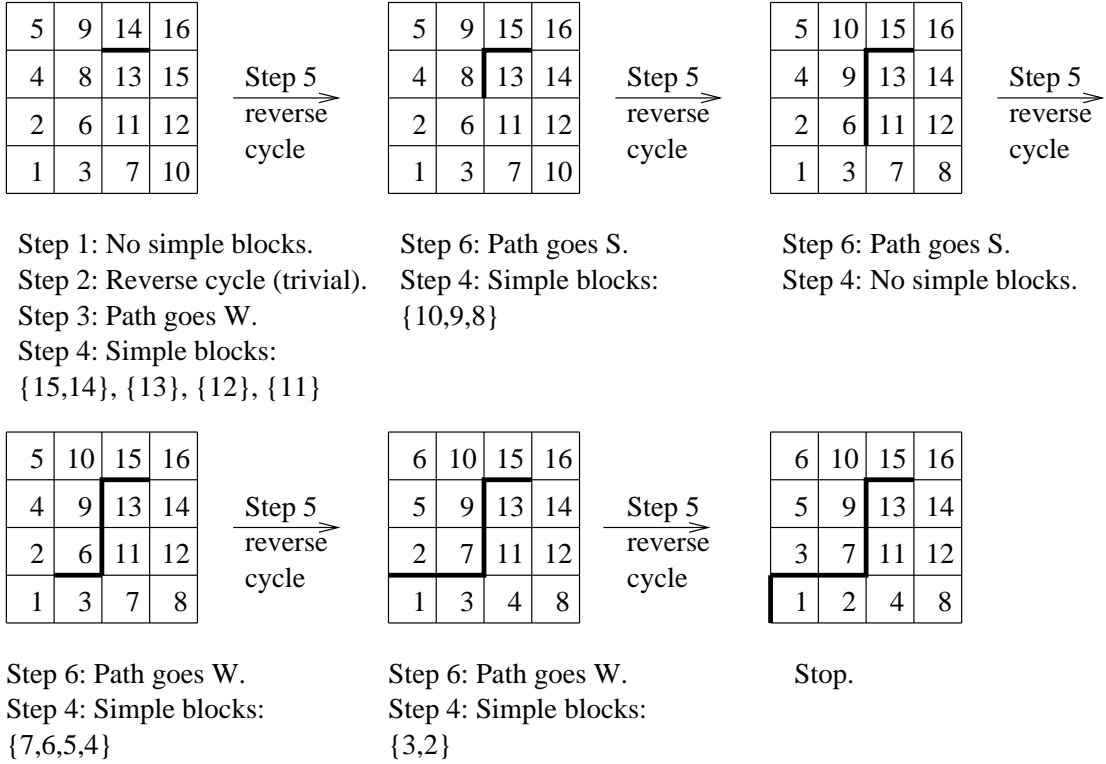}}\end{center}
\caption{An example of the map $\Phi_{n}$}
\end{figure}
\begin{theorem}\label{bij}
For $k=3,\dots,n$, if $T$ is the tableau obtained by applying the above algorithm to the tableau $S$, then $S=\Psi_{k}(T)$.  
\end{theorem}
\begin{proof}
It suffices to prove the result for $k=n$.
First we show that $T\in SYT(\la)$.  As in the proof that $\Psi_{k}:SYT(\la)\to SYT(\la)$,  this follows from the geometry of the simple blocks.  In particular, if the cells with contents $a,\dots,a+m$ form a simple block at any stage of the algorithm, then the cells containing $a,\dots,a+m-1$ are strictly NW, respectively strictly SE, of the cell containing $a+m$ if $n-1\in\Des(S)$, respectively $n-1\notin\Des(S)$.  For example, suppose $n-1\in\Des(S)$ and $c(i,j)=a+m$ and $c(i,j')=a+p$ for some $p$ such that $1\leq p\leq m-1$.  Then $\tilde{\pi}(S,n)$ has a South segment across row $i$.  Thus we must have had, at a previous run through the algorithm, $(i,l)$ in a simple block for some $l$ such that $j'<l\leq j$.   This leads to a contradiction since $c(i,l)\leq a+m$, so $(i,j)$ would have had to have been in a simple block before we constructed the South segment across row $i$.  The rest of the details are left to the reader.

Next we prove that $S=\Psi_{n}(T)$.  It is enough to show that $\tilde{\pi}(S,n)=\pi(T,n)$ because then by Lemma \ref{maj}, the simple blocks in the above algorithm are (up to reordering) the blocks for $\pi(T,n)$ with the exception of the block $\{(1,1)\}$, and the reverse cycling procedure undoes the cycling procedure of the map $\Psi_{n}$.  We give the argument for the case  $n-1\in \Des(S)$; the other case follows by considering the conjugate tableaux.  Suppose we have just completed either Step $2$ or Step $4$ of the algorithm and $\tilde{\pi}(S,n)$ and $\pi(T,n)$ are equal up to this point.  We find ourselves at the lower left-hand corner of a cell $(i,j)$ and we need to show that $\tilde{\pi}(S,n)$ and $\pi(T,n)$ go in the same direction.  This is clear if $i=1$ or $j=1$, so we assume $i,j\neq1$.  Let $p,q$ be the contents of $(i-1,j)$ and $(i,j-1)$, respectively, at this stage.  
Let $c(i-1,j)_{T}=p'$ and $c(i,j-1)_{T}=q'$.  First we consider the case in which $(i-1,j)$ was in a simple block at a previous step in the algorithm.  If $(i,j-1)$ was also in a simple block, then $p=p'$ and $q=q'$ so $\pi(T,n)$ and $\tilde{\pi}(S,n)$ go in the same direction.  If $(i,j-1)$ was not yet in a simple block, then $q<p'=p$ and the cell with content $q'$ will be in a simple block with $(i,j-1)$ at a later stage in the algorithm, so also $q'<p'$.  Thus both paths go South.  Now suppose  $(i-1,j)$ has not yet appeared in a simple block.  Then $c(i-1,j-1)=p-m$ for some $m$ such that the cells with contents $p-1,\dots,p-m+1$ are all above the portion $\tilde{\pi}(S,n)$ constructed so far.  If $(i,j-1)$ was not yet in a simple block,
then $q<p$, but also $q>p-m$, so $q=p-l$ for some $l$ such that $1\leq l\leq m-1$.  By the geometry of the simple blocks, the cell $(i-1,j-1)$ cannot be in the block of $(i-1,j)$, so $p'=p-m+1<q'=p-l+1$.
If $(i,j-1)$ was already in a simple block then $p<q=q'$, and the cell with content $p'$ will be in a simple block with $(i-1,j)$ at a later stage, so also $p'<q'$.  Thus, in both cases, $\pi(T,n)$ and $\tilde{\pi}(S,n)$ go West.  
\end{proof}
Thus $\Psi$ is a bijection with inverse $\Phi=\Phi_{n}\circ\cdots\circ\Phi_{3}$, and it remains only to prove the second statement of Theorem \ref{main}.
\begin{proof}[Proof of Theorem \ref{main}]  Let $T\in SYT(\la)$.  First we show that if $(1,1)$ is below $\pi(T,n)$, then 
$\Inv(T)=\Inv(\Psi_{n}(T)\backslash\{n\})+n-1$, where $\Psi_{n}(T)\backslash\{n\}$ is the standard Young tableau obtained from $\Psi_{n}(T)$ by deleting the cell with content $n$.  We begin by observing that  the set of inversion paths for $\Psi_{n}(T)\backslash\{n\}$ equals the set of inversion paths for $T$ minus the inversion path of the cell containing $n$.  Let $(i,j),(k,l)\in T$ such that $c(k,l)_{T}<c(i,j)_{T}<n$.  First we consider the case in which $(i,j)$ and $(k,l)$ are in different blocks for $\pi(T,n)$.  Because the cycling procedure preserves the relative ordering of the contents of cells within different blocks and the inversion paths of $(i,j)$  
for $T$ and $\Psi_{n}(T)$ are the same, $(i,j)$ and $(k,l)$ form an inversion pair in $T$ if and only if they form an inversion pair in $\Psi_{n}(T)$.  Now suppose $(i,j)$ and $(k,l)$ are in the same block for $\pi(T,n)$.  Say the contents of the cells in this block are $a,\dots,a+m$.  If $(i,j)$ and $(k,l)$ are both above $\pi(T,n)$,  then 
$c(i,j)_{\Psi_{n}(T)}=c(i,j)_{T}-1$ and $c(k,l)_{\Psi_{n}(T)}=c(k,l)_{T}-1$.  Thus, $c(k,l)_{\Psi_{n}(T)}<c(i,j)_{\Psi_{n}(T)}$.   Again because the inversion paths of $(k,l)$ for both tableaux are the same, $(i,j)$ and $(k,l)$ form an inversion pair in $T$ if and only if they form an inversion pair in $\Psi_{n}(T)$.  Finally, suppose $(i,j)$ is above $\pi(T,n)$ and $(k,l)$ is below $\pi(T,n)$. So $c(i,j)_{T}=a+p$ and $c(k,l)_{T}=a$ where $1\leq p\leq m$.  By Lemma \ref{geom}, $(k,l)$ is strictly SE of $(i,j)$, so $(i,j)$ and $(k,l)$ form an inversion pair in $T$.  But $c(i,j)_{\Psi_{n}(T)}=a+p-1$ and  $c(k,l)_{\Psi_{n}(T)}=a+m$, so $(i,j)$ and $(k,l)$ do not form an inversion pair in $\Psi_{n}(T)$.  Thus the total number of inversion pairs in $\Psi_{n}(T)\backslash\{n\}$ equals the number of inversion pairs in $T$ minus one inversion pair for every cell above $\pi(T,n)$ minus one inversion pair for every cell below $\pi(T,n)$, since all of these cells formed inversion pairs in $T$ with the cell with content $n$.  Hence, $\Inv(T)=\Inv(\Psi_{n}(T)\backslash\{n\})+n-1$.

Next we show that If $(1,1)$ is above $\pi(T,n)$, then $\Inv(T)=\Inv(\Psi_{n}(T)\backslash\{n\})$.
Observe that two cells in $T$ form an inversion pair if and only if they do not form an inversion pair in $T'$, the conjugate tableau.  Thus $\Inv(T)={n\choose2}-\Inv(T')$.  Similarly, $\Inv(\Psi_{n}(T)\backslash\{n\})={{n-1}\choose 2}-\Inv((\Psi_{n}(T)\backslash\{n\})')$.  But $(1,1)$ is below $\pi(T',n)$, so
$$\Inv(\Psi_{n}(T)\backslash\{n\})={{n-1}\choose 2}-(\Inv(T')-(n-1))={n\choose2}-\Inv(T')=\Inv(T).$$ 

Thus, by Lemma \ref{maj}
$$\Inv(T)=\Inv(\Psi_{n}(T)\backslash\{n\})+(n-1)\chi\left(n-1\in \Des(\Psi_{n}(T))\right),$$
where for any logical statement $L$, $\chi(L)=1$ if $L$ is true and $\chi(L)=0$ if $L$ is false.
But $\Psi_{k}$ acts as the identity on all cells with content $\geq k$, so
$$\Inv(T)=\Inv(\Psi_{n}(T)\backslash\{n\})+(n-1)\chi\left(n-1\in \Des(\Psi(T))\right).$$
Iterating this argument, we see that
$$\Inv(T)=\sum_{j=1}^{n-1}j\chi\left(j\in \Des(\Psi(T))\right)=\maj(\Psi(T)).$$
\end{proof}

In Figure $6$, we give an example of a tableau $S$ and its image under the map $\Psi$.  We have $\maj(\Psi(S))=1+4+5+7=17$.  In Figure $7$, we show the set of inversion paths for $S$.  If we label cells by their contents, then inversion pairs for $S$ are $(8,6),(8,5),(8,4),(8,3),(8,2),(8,1),(7,6),(7,5),(7,4),(7,3),(7,2),(7,1),(6,5),(4,2),$\newline $(4,1),(3,2),(3,1).$  So $\Inv(S)=17$.

\begin{figure}\label{Phi}
\begin{center}
{\includegraphics{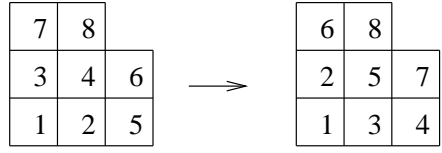}}\end{center}
\caption{On the left: a tableau $S$.  On the right: $\Psi(S)$.}
\end{figure}

\begin{figure}\label{allpaths}
\begin{center}
{\includegraphics{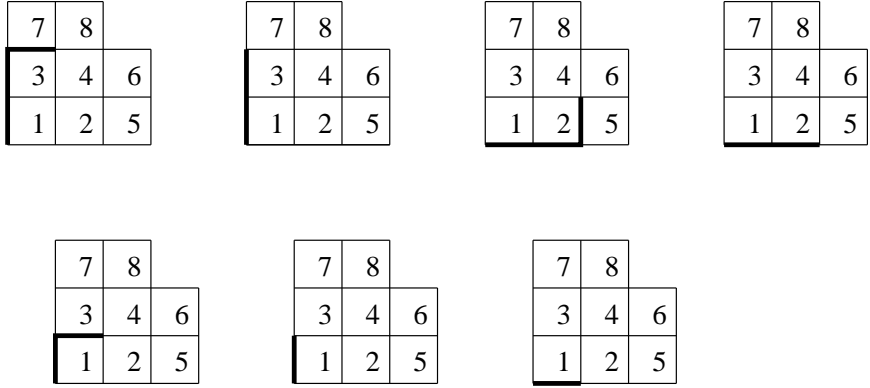}}\end{center}
\caption{The set of inversion paths for $S$.}
\end{figure}

\begin{corollary}
\label{majinv}
The statistics $\maj$ and $\Inv$ have equal distributions over $SYT(\la)$.
\end{corollary} 

\noindent{\textbf{Remark:}} 
Part of our motivation for defining an inversion statistic on $SYT$
is based on recent work of 
Haglund, Haiman and Loehr \cite{HHL}, where
it is shown that the coefficient of a monomial symmetric function in 
the Macdonald polynomial ${\tilde H}_{\mu}[X;q,t]$ can be expressed in terms of 
descents and inversion pairs for fillings of the Ferrers shape of $\mu$ by
positive integers.   One could hope that a similar phenomenon would apply to the
coefficients ${\tilde K}_{\lambda ,\mu}(q,t)$ in the expansion of 
${\tilde H}_{\mu}$ into Schur functions (Macdonald showed that
${\tilde K}_{\lambda ,\mu}(1,1)$ equals the number of $SYT(\lambda)$, and posed the question \cite[p. 356]{Mac} of whether 
${\tilde K}_{\lambda ,\mu}(q,t)$ could be written in the form
\begin{align}
\label{Mac}
{\tilde K}_{\lambda ,\mu}(q,t) = \sum_{T \in SYT(\lambda)}
t^{\text{tstat}(T,\mu)}
q^{\text{qstat}(T,\mu)}
\end{align}
for some tableau statistics $\text{tstat},\text{qstat}$).  
In \cite[Conjecture 3]{Hag04} specific values for $\text{tstat},\text{qstat}$ are
conjectured for all $\mu$ with at most three columns.  The value of 
$\text{tstat}$ is described in terms of descents of $T$, while the value of $\text{qstat}$ is
described in terms of certain ``inversion pairs" of $T$.  (The special case of this 
conjecture where $\mu$ has at most two columns is proved in
\cite[Proposition 9.2]{HHL}).  
Although we have as yet been unable to extend this 
conjecture to general $\mu$, Corollary \ref{majinv} gives a way of expressing 
${\tilde K}_{\lambda, (n)}(q,t)$ in this form, since it is known that
\begin{align}
\label{Mac2}
{\tilde K}_{\lambda ,(n)}(q,t) = \sum_{T \in SYT(\lambda)}
q^{\text{maj}(T)},
\end{align}
which by Corollary \ref{majinv} also equals 
\begin{align}
{\tilde K}_{\lambda ,(n)}(q,t) = \sum_{T \in SYT(\lambda)}
q^{\text{Inv}(T)}.
\end{align}

\section{Extensions}
\label{extensions}

Although we have assumed throughout this article
that $\lambda$ is a partition shape, all our
results apply just as easily to skew shapes.
For a given SYT $T$ of skew shape $\lambda \vdash n$, note  
that our map $\Phi $ (and its inverse) both fix the largest element $n$ in the tableau.
Thus our inversion statistic is equidistributed with $\text{maj}$ over the set of
$SYT(\lambda)$ with $n$ occurring in some fixed corner square.  Let
$$
\text{comaj}(T) = \sum_{i \in \text{Des}(T)} n-i.
$$
It is known that $\text{comaj}$ and $\text{maj}$ have the same distribution over
$SYT$ of skew shape \cite[Chapter 7]{Sta}.
We briefly describe a parallel version of our Foata map which is linked to 
$\text{comaj}$.  If square $(i,j)$ contains $1$ in $T$, 
start at the NE corner of $(i,j)$ and draw a 
unit line segment $N$ if $c(i+1,j)>c(i,j+1)$, otherwise draw a unit line segment $E$,
and now iterate, eventually ending up at the upper border of $T$.
Next break up the numbers $2$ through $n$ into maximal blocks, as before, 
consisting of an integer 
$k$ on the same side of the inversion path as $n$, and a consecutive sequence of integers
$k-1,k-2,\ldots, k-a$ on the other side of the path, then cycling by moving $k-a$ to the square containing $k$, $k$ to the square containing $k-1$, etc..  
Thus each square ends up with an inversion path which starts at its $NE$ corner and 
travels $NE$ to the border, and the number of inversion pairs are equidistributed with
$\text{comaj}$ over the set of all $SYT(\lambda)$ with $1$ occurring in some fixed corner square.

We now show that the special case of our map $\Phi$ 
when $\lambda$ is a disjoint union of
squares, i.e. $\lambda = (n,n-1,\ldots ,1)/(n-1,\ldots ,1)$, and $T$ is the
tableau obtained by filling these squares with the numbers $\sigma _1, \ldots ,
\sigma _n$ (see Figure $8$), is essentially
the same as Foata's original map $\phi$ applied to $\sigma ^{-1}$.  We identify $T$ with $\sigma$ and use the notation $\Phi(\sigma)$.  To find $\Phi_{k}(\sigma)$, partition the numbers $k-1,\dots,1$ into blocks $\{a,a-1,\dots,a-j\}$ of maximal length such that $a$ is on the same side of $k$ as $k-1$ and $a-1,\dots,a-j$ are on the other side of $k$.  Then reverse cycle, e.g. $\Phi_{4}(346251)=146352$.
Let $\sigma^{-1}=\omega$, so $\omega_{i}$ is the position of $i$ in $\sigma$.  We now show that
$\phi(\omega)=(\Phi(\sigma))^{-1}$.  The main idea is to prove that for any $k,\, 3\leq k\leq n$, 
$\phi^{(k)}\omega_{k+1}\dots\omega_{n}=(\Phi_{k}\circ\cdots\circ\Phi_{3}(\sigma))^{-1}$, where $\phi^{(k)}=\phi^{(k)}(\omega)$.  This is straightforward when $k=3$; $\Phi_{3}(\sigma)=\sigma$ unless $2$ and $1$ are on opposite sides of $3$, in which case $\Phi_{3}$ interchanges their positions and leaves all other numbers fixed, whereas $\phi^{(3)}=\omega_{1}\omega_{2}\omega_{3}$ only if $\omega_{1}$ and $\omega_{2}$ are both less than or greater than $\omega_{3}$, otherwise $\phi^{(3)}=\omega_{2}\omega_{1}\omega_{3}$.  Now assume that the claim is true for some $k$ such that $3\leq k\leq n-1$.  Write $\Phi_{k}\circ\cdots\circ\Phi_{3}(\sigma)=\beta$ and $\phi^{(k)}=\alpha_{1}\dots\alpha_{k}$.  Then $\{a,a-1,\dots,a-j\}$ is a block in $\beta$ if and only if $\alpha_{a-j}\dots\alpha_{a-1}\alpha_{a}$ is a block in $\alpha_{1}\dots\alpha_{k}\omega_{k+1}$.  The reverse cycling procedure sends $a-j$ to the position of $a$ and $a,\dots,a-j+1$ to the positions of $a-1,\dots,a-j$, respectively, which corresponds to the Foata algorithm sending the block $\alpha_{a-j}\dots\alpha_{a-1}\alpha_{a}$ to $\alpha_{a}\alpha_{a-j}\dots\alpha_{a-1}$.  Thus $\phi^{(k+1)}\omega_{k+1}\dots\omega_{n}=(\Phi_{k+1}(\beta))^{-1}.$  For the example in Figure $8$, the sequence $\sigma,\Phi_{3}(\sigma),\Phi_{4}\circ\Phi_{3}(\sigma),\Phi_{5}\circ\Phi_{4}\circ\Phi_{3}(\sigma),\Phi(\sigma)$ is
$$346251,346251,146352,146253,256143$$ whereas the sequence $\omega,\phi^{(3)}\omega_{4}\omega_{5}\omega_{6},\phi^{(4)}\omega_{5}\omega_{6},\phi^{(5)}\omega_{6},\phi(\omega)$ is
$$641253,641253,164253,146253,416523.$$

\begin{figure}\label{skew}
\begin{center}
{\includegraphics{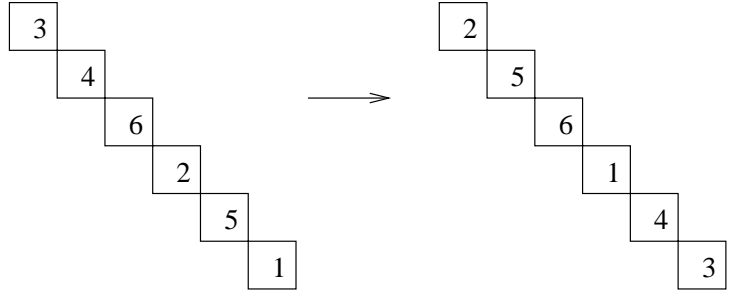}}\end{center}
\caption{On the left: a tableau $T$ of shape $(6,5,4,3,2,1)/(5,4,3,2,1)$ whose filling corresponds to the permutation $\sigma=346251$.  On the right: $\Phi(T)$.}
\end{figure}

\section{Concluding Remarks}
\label{concluding remarks}

The authors would like to thank the referees for helpful suggestions on the exposition, as well as pointing out some interesting directions for future research.  
For example, MacMahon's result 
on the equidistribution of $\text{maj}$ and $\text{inv}$ holds for multiset 
permutations, and Foata showed how his map proves this more general fact
bijectively.  Perhaps our results on $\text{SYT}$ also have versions for 
$\text{SSYT}$.  We remark that the obvious thing to try, which is to start with a 
$\text{SSYT}$, standardize in the usual way, apply the $\Phi$ map, then ``unstandardize", doesn't quite work, as this 
process will in general produce fillings which are not column strict, and hence not
$\text{SSYT}$.  One thing that does seem to have a
natural analogue in our setting is that of the ``inversion code" of a permutation, 
which is a sequence
$x_1 x_2\ldots x_{n}$ with $0\le x_i \le i-1$ for $1\le i \le n$, where
$x_i$ equals the number of inversion pairs of the permutation of the
form $(i,j)$ with $i>j$.
Our tableaux inversion
statistic is exactly the number of similar such pairs, as in the example just 
above Corollary \ref{majinv}.  One could also hope to connect our ideas with other inversion statistics which occur in the theory of symmetric functions and Macdonald polynomials, such as the
inversion statistic for tuples of $\text{SSYT}$ which can be used to
describe LLT polynomials (see  
\cite{HHL}).

\end{document}